\documentclass[12pt]{article}
\usepackage[english]{babel}
\usepackage{amssymb,amsmath,amsthm}
\usepackage{graphicx}
\usepackage{color}
\theoremstyle{plain}
\newtheorem{de}{Definition}
\newtheorem{theor}{Theorem}
\newtheorem{pr}{Proposition}

\newtheorem{example}{Example}

\begin{document}

\begin{center}
{\bfseries Rigidity of some classes of Lie algebras in connection to Leibniz algebras}

\vspace{0.5cm}
\underline{Abdulafeez O. Abdulkareem$^{1}$, Isamiddin S. Rakhimov}$^{2}$, and Sharifah K. Said Husain$^{3}$   \\[0.5cm] 
$^{1,3}${\it Department of Mathematics, Faculty of Science, \\ Universiti Putra Malaysia, 43400 UPM
  Serdang, \\ Selangor Darul Ehsan, Malaysia.} \\[0.5cm]
$^{1,2,3}${\it Laboratory of Cryptography,Analysis and Structure, \\ Institute for Mathematical Research (INSPEM), \\ Universiti Putra Malaysia, 43400 UPM  Serdang, \\ Selangor Darul Ehsan, Malaysia.} \\[0.5cm] 
\end{center}

\begin{abstract}
In this paper we focus on algebraic aspects of contractions of Lie and Leibniz algebras. The rigidity of algebras plays an important role in the study of their varieties. The rigid algebras generate the irreducible components of this variety. We deal with Leibniz algebras which are generalizations of Lie algebras. In Lie algebras case, there are different kind of rigidities (rigidity, absolutely rigidity, geometric rigidity and e.c.t.). We explore the relations of these rigidities with Leibniz algebra rigidity. Necessary conditions for a Lie algebra to be Leibniz rigid are discussed.
\end{abstract}

\section{Introduction}
In 1951, Segal I.E. \cite{155} introduced the notion of contractions of Lie algebras on physical grounds: if two physical theories (like relativistic and classical mechanics) are related by a limiting process, then the associated invariance groups (like the Poincar\'e and Galilean groups) should also be related by some limiting process. If the velocity of light is assumed to go to infinity, relativistic mechanics "transforms" into classical mechanics. This also induces a singular transition from the Poincar\'e algebra to the Galilean one. Another example is a limiting process from quantum mechanics to classical mechanics under $\hbar \to 0,$  that corresponds to the contraction of the Heisenberg algebras to the abelian ones of the same dimensions \cite{41}.\\
There are two approaches to the contraction problems of algebras. The first of them is based on physical considerations that is mainly oriented to applications of contractions. Contractions were used to establish connections between various kinematical groups and to shed a light on their physical meaning. In this way relationships between the conformal and Schrodinger groups was elucidated and various Lie algebras including a relativistic position operator were interrelated. Under dynamical group description of interacting systems, contractions corresponding to the coupling constant going to zero give noninteracting systems. Application of contractions allows to derive interesting results in the special function theory as well. The second consideration is purely algebraic, dealing with abstract algebraic structures. We will deal with this case and focus mainly on algebraic aspects of the contractions.\\
Let $V$  be a vector space of dimension $n$  over an algebraically closed field $K$ ($charK =0$). The set of bilinear maps $V\times V\longrightarrow V$ forms a vector space $Hom(V\otimes V, V)$ of dimension $n^{3}$, which can be considered together with its natural structure of an affine algebraic variety over $K$. It is denoted by $Alg_{n}(K)$. An $n$-dimensional algebra $A$ over $K$ may be regarded as an element $\lambda(A)$ of $Alg_{n}(K)$ via the bilinear mapping $\lambda:A\otimes A\longrightarrow A$ defining a binary algebraic operation on $A$. The linear reductive group $G=GL_{n}(K)$ acts on $Alg_{n}(K)$ by $(g *\lambda)(x,y)=g(\lambda(g^{-1}(x),g^{-1}(y)))$  (“transport of structure”). Two algebras $\lambda_{1}$ and $\lambda_{2}$ are isomorphic if and only if they belong to the same orbit under this action. For given two algebras $\lambda$ and $\mu$ we say that $\lambda$ degenerates to $\mu$, if $\mu$ lies in the Zariski closure of the orbit $\lambda$. We denote this by $\lambda\to \mu$. 
\section{Preliminaries}
\begin{de}
An algebra $L$ over a field $K$ is called a Leibniz algebra if its binary operation $\lambda$ satisfies the following Leibniz identity: 
\begin{equation} \lambda(x,\lambda(y,z))=\lambda(\lambda(x,y),z)-\lambda(\lambda(x,z),y), \ for \ any \ x,y,z\in L.
\end{equation}
\end{de}
The set of all $n$-dimensional Leibniz algebras over a field $K$  will be denoted by $LB_{n}(K)$ . The set $LB_{n}(K)$ can be included in the above mentioned $n^{3}$-dimensional affine space as follows: let $\{e_{1},e_{2},\cdots, e_{n}\}$ be a basis of the Leibniz algebra $L$, then the table of multiplication of $L$ is represented by point $\gamma^{k}_{ij}$ of this affine space as follows: $$\lambda(e_i, e_j) =\sum \limits_{k=1}^n \gamma_{ij}^k e_k, \ \ i,j=1,2,...,n.$$
Thus, the algebra $L$ corresponds to the point $\gamma^{k}_{ij}$ which are called structure constants of $L$. The Leibniz identity gives polynomial relations among $\gamma^{k}_{ij}$. Hence we regard $LB_{n}$ as a subvariety of $K^{n^{3}}$.
\begin{de}
A Leibniz algebra $\lambda$ is said to degenerate to a Leibniz algebra $\mu$ if $\mu$ is represented by a structure which lies in the Zariski closure of the $GL_{n}(K)$-orbit of the structure which represents $\lambda$.
\end{de}
In this case entire orbit $O(\mu)$ lies in the closure of $O(\lambda)$. We denote this, as has been mentioned above, by $\lambda\to \mu$ i.e., $\mu\in \overline{O(\lambda)}$. Degeneration is transitive, that is if $\lambda\to \mu$ and $\mu\to \nu$ then $\lambda\to \nu$. There are algebras whose orbits are open in $LB_{n}(K)$. These algebras are called rigid. In that case, the corresponding algebra does not admit any non-trivial deformation. The orbits of the rigid algebras give irreducible components of the variety $LB_{n}(K)$. However, in general, not every irreducible component is generated by rigid algebras. Flanigan F.J. \cite{54} has shown that there exists a component in $Alg_{3}$ which consists of union of infinitely many orbits of non isomorphic algebras having the same dimension. If $L$ is a rigid algebra in $LB_{n}(K)$ then, there exists an irreducible component $\mathfrak{C}$ of $LB_{n}(K)$ such that $O(L)\cap \mathfrak{C}$ is non empty open subset of $\mathfrak{C}$. The closure of $O(L)$ is contained in $\mathfrak{C}$. Then, the dimension of the irreducible component $\mathfrak{C}$ of $LB_{n}(K)$ is given by $$dim\mathfrak{C}=n^{2}-dimAutL.$$
The problem of counting the number of the irreducible components and the number of open orbits of the variety of algebras is of interest (see \cite{80}, \cite{104}). Here is a proposition from \cite{26}. In Lie algebras case, counting the number of irreducible components and the number of open orbits.
\begin{pr}
Let $r(n)$ and $s(n)$ be the number of irreducible components and the number of open orbits, respectively, in the variety $L_{n}(\mathbb{C})$ of $n$-dimensional complex Lie algebras. Then one has $$(r(1),r(2),\cdots,r(7))=(1,1,2,4,7,17,49)$$
and 
$$(s(1),s(2),\cdots,s(7))=(1,1,1,2,3,6,14)$$
\end{pr}
Hence to find the rigid orbits of variety of algebras is of great interest. By the Noetherian consideration they are finite number.\\
From now on all algebras considered are supposed to be over the field of complex numbers $\mathbb{C}$ unless otherwise specified. We make use of a few useful facts from the algebraic groups theory, concerning the degenerations. The first of them is on constructive subsets of algebraic varieties over $\mathbb{C}$, the closures of which relative to Euclidean and Zariski topologies coincide. Since $GL_{n}$-orbits are constructive sets, the usual Euclidean topology on $\mathbb{C}^{n^{3}}$ leads to the same degenerations as does the Zariski topology. Now we may express the concept of degeneration in a slightly different way, that is the following condition will imply that $\lambda\to \mu$: $$\exists g_{t}\in GL_{n}(\mathbb{C}(t)) \ such \ that \lim\limits_{t\to 0}g_{t}* \lambda=\mu$$
where $\mathbb{C}(t)$ is the field of fractions of the polynomial ring $\mathbb{C}[t]$
The second fact concerns the closure of  $GL_{n}(\mathbb{C})$-orbits stating that the boundary of each orbit is a union of finitely many orbits with dimensions strictly less than dimension of the given orbit. It follows that each irreducible component of the variety, on which algebraic group acts, contains only one open orbit that has a maximal dimension. It is obvious that in the content of variety of algebras the representatives of such orbits are rigid.
It is an interesting but difficult problem to determine the number of irreducible components of an algebraic variety. But if one is interested in finding the dimension of the algebraic variety then degeneration approach is also helpful. In this case no need to find all the degenerations, just to find the so-called rigid algebras. The closure of the rigid algebra gives a component of the variety. In order to find the dimension of the variety it is sufficient to find a rigid algebra having a maximal orbit dimension. 
Let $A$ be an $n$-dimensional algebra (underlying vector space is denoted by $V$) with the binary operation   Consider a continuous function $g_{t}:(0,1]\longrightarrow GL(V)$. In other words, $g_{t}$ is a nonsingular linear operator on $V$ for all $t\in (0,1]$. Define parameterized family of new, isomorphic to $A$, algebra structures on $V$ via the old binary operation $\lambda$ as follows:
$$\lambda_{t}(x,y)=(g_{t}*\lambda)(x,y)=g_{t}^{-1}\lambda(g_{t}(x),g_{t}(y))$$  
\begin{de}
 If the limit $\lim\limits_{t\to +0}\lambda_{t}=\lambda_{0}$ exists for all $x, y\in V$ then, the algebraic structure $\lambda_{0}$ defined by this way on $V$ is said to be a contraction of the algebra Obviously, the contractions can be considered in basis level, i.e, let $\{e_{1},e_{2},\cdots,e_{n}\}$ be a basis of an $n$-dimensional algebra $A$. If the limit $\lim\limits_{t\to+0}\lambda_{t}(e_{i},e_{j})=\lambda_{0}(e_{i},e_{j})$ exists then the algebra $(V,\lambda_{o})$ is a contraction of  
$A$.
\end{de}
\begin{de}
A contraction from an algebra $A$ to algebra $A_{0}$ is said to be trivial if $A_{0}$ is abelian and improper if $A_{0}$ is isomorphic to $A$ 
\end{de}
Note that both the trivial and the improper contractions always exist. Here is an example of the trivial and the improper contractions.
\begin{example}
 Let $A=(V,\lambda)$ be an $n$-dimensional algebra. If we take $g_{t}=diag(t,t,\cdots,t)$ then $g_{t}* \lambda$ is abelian and at $g_{t}=diag(1,1,\cdots,1)$  we get $g_{t}* \lambda=A$. 
\end{example} 
We equally use the notions of rigidity with respect to degeneration and deformation synonymously since the rigidity with respect to degeneration follows from the rigidity with respect deformation, but the converse does not always hold. In the case of Lie algebras over a field characteristic zero, the two definitions of rigidity are equivalent \cite{51}, \cite{98}. By using the results of \cite{51}, this can be extended to the Leibniz algebras case.\\
Algebraic deformation theory was introduced, for associative algebras, in 1963 by Gerstenhaber et. al. \cite{61} and it was extended to Lie algebras by Nijenhuis and Richardson \cite{105}. Following Gerstenhaber, Nijenhuis and Richardson, many authors have published papers on some aspect of the deformations of a given type of algebraic structures (associative, commutative, Lie, Jordan,...).\\
One fundamental concept in deformation theory is the notion of rigidity. Loosely speaking, an algebra $A$ is rigid if an arbitrary "infinitesimal" deformation of $A$ produces an algebra $A$ isomorphic to $A$. Focusing on Lie algebras, a geometric formulation of rigidity is given as follows.\\
Let $V$ be a complex vector space of dimension $n$. Let $L_{n}$ denotes the set of all Lie algebra structures on $V$. The set $L_{n}$ has a natural structure of an affine variety. There is an action of the group $GL_{n}$, of all automorphisms of $V$, on $L_{n}$. Under this action each orbit represents an isomorphism class of Lie algebra structures. A Lie algebra $L$ is rigid if its orbit under the action of $GL_{n}$ is Zariski open in $L_{n}$. A theorem of Nijenhuis and Richardson \cite{105} shows that the vanishing of the second cohomology of a Lie algebra, with coefficients in the adjoint representation, is a sufficient condition for a Lie algebra to be rigid. Lie algebras are part of a bigger class of algebras called Leibniz algebras. The latter were introduced by Loday \cite{88}. Then Balavoine \cite{14} published a paper on the deformations of Leibniz algebras. In particular, using the theory of Fox \cite{56}, he obtained a result analogous to that of Nijenhuis and Richardson, i.e., the vanishing of the second Leibniz cohomology is a sufficient condition for the rigidity of a Leibniz algebra. Since a Lie algebra is also Leibniz algebra, we need to study the rigidity of Lie algebras as Lie algebras and as Leibniz algebras.
\section{On Rigidity of Lie Algebras}
Let $V$ be a vector space over a field $K$. A Lie algebra multiplication on $V$ is an element $\mu$ of $Hom(V\otimes V, V)$ satisfying:
\begin{enumerate}
\item[(i)] $\mu(x,y)=-\mu(y,x)$ for all $x,y \in V$
\item[(ii)] $\mu(x,\mu(y,z))+\mu(y, \mu(z,x))+\mu(z,\mu(x,y))=0$ for all $x,y,z\in V$.
\end{enumerate} 

The pair $L=(V,\mu)$ is called a Lie algebra over $K$. The condition (ii) is called the Jacobi identity. 
\begin{de}
Let $L=(V,\mu)$ be a Lie algebra. A Lie module is a vector space $M$ together with an action   satisfying $$\mu(x,y)\cdot m=x\cdot (y\cdot m)-y\cdot(x\cdot m)$$ for $x,y\in V$ and $m\in M$.
\end{de}
\subsection{(Co)homology of Lie Algebras}
Let $L=(V,\mu)$ be a Lie algebra and $M$ be a Lie module. Let $C^{n}(L,M)=Hom(\Lambda^{n}V,M)$, for $n\geq 0$ define $d:C^{n}(L,M)\longrightarrow C^{n+1}(L,M)$ by 
\begin{align*}
(df)(x_{1},x_{2},\cdots,x_{n+1})&=\sum_{i=1}^{n+1}(-1)^{i+1}x_{i}\cdot f(x_{1},\cdots,x_{i-1},\hat{x_{i}},x_{i+1},\cdots,x_{n+1})+\\&\sum_{1\leq i,j\leq n+1}(-1)^{i+j}f(\mu(x_{i},x_{j}),x_{1},\cdots,x_{i-1},\hat{x_{i}},x_{i+1},\cdots,x_{j-1},\hat{x_{j}},x_{j+1},\cdots,x_{n+1}),
\end{align*}
where $\hat{x_{i}}$ and $\hat{x_{j}}$ indicate omitted terms and ``$\cdot$" is the module action. The map $d$ must have the property that $d\circ d=0$, which is equivalent to antisymmetricity of $\Lambda,\mu$ and the Jacobi identity.\\
Let $B^{n}(L,M)=\{f\in C^{n}(L,M)|f=dg \ for \ some \ g\in C^{n-1}(L,M)\}$ and $Z^{n}(L,M)=\{f\in C^{n}(L,M)|df=0\}$. $B^{n}(L,M)$ is called the space of all coboundaries while $Z^{n}(L,M)$ is the space of all $n$-cocycles. The fact that $d\circ d=0$ implies that $B^{n}(L,M)$ is a subspace of $Z^{n}(L,M)$. The quotient space $H^{n}(L,M)=Z^{n}(L,M)/B^{n}(L,M)$ is called the $nth$-cohomology space of the Lie algebra $L$ with coefficients in the Lie module $M$.\\
In the particular case when $M=L$ under the adjoint action, i.e.,$L\times L\longrightarrow L(x,y)\mapsto \mu(x,y)$, the coboundary operator $d$ is written as 
\begin{align*}
(df)(x_{1},x_{2},\cdots,x_{n+1})&=\sum_{i=1}^{n+1}(-1)^{i+1}\mu(x_{i}, f(x_{1},\cdots,x_{i-1},\hat{x_{i}},x_{i+1},\cdots,x_{n+1}))+\\&\sum_{1\leq i,j\leq n+1}(-1)^{i+j}f(\mu(x_{i},x_{j}),x_{1},\cdots,x_{i-1},\hat{x_{i}},x_{i+1},\cdots,x_{j-1},\hat{x_{j}},x_{j+1},\cdots,x_{n+1}),
\end{align*}

Note that $d:C^{1}(L,L)\longrightarrow C^{2}(L,L)$ is given by; $$(df)(x,y)=\mu(x,f(y))-\mu(y,f(x))-f(\mu(x,y))$$ and $d:C^{2}(L,L)\longrightarrow C^{3}(L,L)$ is given by; $$(df)(x,y,z)=\mu(x,f(y,z))-\mu(y,f(x,z))+\mu(z,f(x,y))-f(\mu(x,y),z)+f(\mu(x,z),y)-f(\mu(y,z),x).$$

Let $L_{n}$ be the set of all Lie algebra structures on an $n$-dimensional complex vector space $V$. The set $L_{n}$ has the structure of a complex affine algebraic variety (in general, it is non reduced). The general linear group $GL_{n}$ acts on the set $L_{n}$, and the action is given by
$$(g*\mu)(x,y)=g(\mu(g^{-1}(x),g^{-1}(y))) \ for \ g\in GL_{n} \ and \ \mu\in L_{n}.$$
The orbit $O(\mu)=\{g*\mu|g\in GL_{n}\}$ of $GL_{n}$ on $L_{n}$, for a fixed Lie algebra $L=(V,\mu),$ represents the set of all Lie algebra multiplications such that the Lie algebras $L=(V,\mu)$ and $L'=(V,g*\mu)$ are isomorphic. The closure of $O(\mu)$  with respect to Zariski topology is denoted by $\overline{O(\mu)}$. For given two algebras $\lambda$ and $\mu$ we say that $\lambda$  degenerates to $\mu$ if $\mu$ lies in the Zariski closure of the orbit $\lambda$. The degeneration is said to be trivial if the algebras $\lambda$ and $\mu$ are isomorphic. 
\begin{de}
A Lie algebra $L=(V,\mu)$ is called rigid if its orbit $O(\mu)$ is open in the Zariski topology.
\end{de} 
\begin{theor}(Nijenhuis-Richardson)\label{NR}
 Let $L$ be a Lie algebra such that $H^{2}(L,L)=0$ then $L$ is rigid. 
\end{theor}
The following example from \cite{144} shows that the converse of Theorem \ref{NR} is false.
\begin{example}
Consider three-dimensional semi-simple Lie algebra $Sl_{s}(\mathbb{C})$ over $\mathbb{C}$  and for each positive integer $n$, let
$$\rho_{n}:Sl_{2}(\mathbb{C})\longrightarrow gl(\mathbb{C}^{2n+1})$$ 
be the irreducible representation of weight $n$. 
Regarding $\mathbb{C}^{2n+1}$ as an abelian Lie algebra, let $L_{n}=Sl_{2}+\rho_{n}\mathbb{C}^{2n+1}$ be the semi-direct product of the Lie algebras $Sl_{2}$ and $\mathbb{C}^{2n+1}$ with respect to the representation $\rho_{n}$. It is shown in \cite{144} that $L_{n}$ is rigid for $n\neq 1,2,3,5$ and that $H^{2}(L,L)\neq 0$ if $n$ is odd. This yields a family $L_{7},L_{9},L_{11},\cdots$ of rigid Lie algebras with $H^{2}(L,L)\neq 0$.
\end{example} 
\begin{de}
A Lie algebra $L$ is called absolutely rigid if $H^{2}(L,L)=0$. 
\end{de}
\begin{example}
Semi-simple Lie algebras are absolutely rigid.  
\end{example}
\begin{example}
The general Linear algebra $gl_{n}$ is absolutely rigid.
\end{example}
\begin{example} 
The non abelian two dimensional Lie algebra with the table of multiplication $e_{1}e_{2}=e_{2}$ is absolutely rigid. 
\end{example}
\begin{example}
 Any algebra isomorphic to a Borel subalgebra of a semi-simple Lie algebra is absolutely rigid.
\end{example}

\section{Leibniz Rigidity of Lie algebras}
\subsection{(Co)homology of Leibniz Algebras}
\begin{de}
 Let $(L,[\cdot,\cdot])$ be a Leibniz algebra. Let $M$ be a vector space equipped with two actions (left and right) of $L$, i.e., $[\cdot,\cdot]:L\otimes M\longrightarrow M$   and $[\cdot,\cdot]:M\otimes L\longrightarrow M$ satisfying:
 \begin{enumerate}
 \item[(MLL)] $[m,[x,y]]=[[m,x],y]-[[m,y],x]$
 \item[(LML)] $[x,[m,y]]=[[x,m],y]-[[x,y],m]$
 \item[(LLM)] $[x,[y,m]]=[[x,y],m]-[[x,m],y]$
 \end{enumerate}
 for any $m\in M$ and $x,y\in L.$ Then $M$ is called a Leibniz module over $L$.
\end{de}
The algebra $L$ itself is a Leibniz module under the action $L\otimes L\longrightarrow L$ by $(x,y)\mapsto [x,y]$. In this case the two actions are given by the Leibniz multiplication.
The underlying field $K$ (in our case it is $\mathbb{C}$) also is a Leibniz module under the actions $K\otimes L$ by $(k,x)\mapsto [k,x]=0$ and $L\otimes K$ by $(x,k)\mapsto [x,k]=0$.\\
Let $(L,[\cdot,\cdot])$ be a Leibniz algebra and $M$ be a Leibniz module over $L$. Denote $E^{n}(L,M):=Hom(L^{\otimes n},M)$, $n\geq 0$. The Loday coboundary map $d:E^{n}(L,M)\longrightarrow E^{n+1}(L,M)$ is a homomorphism defined as follows:
\begin{align*}
(df)(x_{1},x_{2},\cdots,x_{n+1})&=[x_{1},f(x_{2},\cdots,x_{n+1})]+\sum_{i=2}^{n+1}(-1)^{i}[f(x_{1},\cdots,\hat{x_{i}},\cdots,x_{n+1}),x_{i}]+\\&\sum_{1\leq i,j\leq n+1}(-1)^{j+1}f(x_{1},\cdots,x_{i-1},[x_{i},x_{j}],x_{i+1},\cdots,\hat{x_{j}},x_{j+1},\cdots,x_{n+1}),
\end{align*}
The condition $d\circ d=0$, $n\geq 0$ is precisely the fact that the Leibniz bracket satisfies the Leibniz identity \cite{39}. Then $(E^{*}(L,M),d)$ is a cochain complex, whose cohomology is called the cohomology of the Leibniz algebra $L$ with coefficients in the representation $M:HL^{*}(L,M):=H^{*}((C^{*}(L,M),d))$. Similarly, we have a chain complex $(E_{*}(L,M),d)$, whose homology is called the homology of the Leibniz algebra $L$ with coefficients in the representation $M:HL_{*}(L,M):=H_{*}((C_{*}(L,M),d))$.
Note that for a Leibniz algebra $L=(V,\mu)$ the operator $d:E^{1}(L,L)\longrightarrow E^{2}(L,L)$ is given by
$$(df)(x,y)=\mu(x,f(y))-\mu(y,f(x))-f(\mu(x,y)).$$ The operator $d:E^{2}(L,L)\longrightarrow E^{3}(L,L)$ is given by; $$(df)(x,y,z)=\mu(x,f(y,z))-\mu(y,f(x,z))+\mu(z,f(x,y))-f(\mu(x,y),z)+f(\mu(x,z),y)-f(\mu(y,z),x).$$
 
For the further references we remind a few known facts about Leibniz cohomology and Leibniz homology, which can be found in \cite{39} and \cite{107}.
\begin{pr}
\begin{enumerate}
\item[(i)] If $L=(V,\mu)$ is an abelian Leibniz algebra, then 
$$HL_{n}(L,M)=M\otimes V^{n}$$
\item[(ii)] For a Leibniz algebra $L$ and a module $M$ over $L$ one has 
$$HL_{n}(L,M)=Hom(HL_{n}(L,M),M)$$
\item[(iii)] If   is a Lie algebra over a field, then  $$HL_{n}(L,L)=HL_{n+1}(L,K)$$ \hfill (\cite{107},p.708)
\item[(iv)] If $L$ is the non abelian two dimensional Lie algebra over a field $K$ and $K$ contains the rational numbers, then
$$HL_{n}(L,K)=K \ for \ all \ n>0$$
\item[(v)] If $L$ is a semi-simple Lie algebra, then 
$$HL_{n}(L,L)=HL^{n}(L,L)=0 \ for \ all \ n>1$$
\item[(vi)] If $L=gl_{n}(K)$ is the general linear algebra, then $$HL_{n}(L,K)=K \ for all \ n>0$$
\end{enumerate}
\end{pr} 
Here is a result from \cite{108} on Leibniz rigidity of Lie algebras. 
\begin{theor}
 Let $L=(V,\mu)$ be a complex Lie algebra. Then
\begin{enumerate}
\item[(a)] $H^{2}(L,L)\subset HL^{2}(L,L)$ 
\item[(b)] The Lie algebra $L$ is Leibniz rigid if and only if $L$ is a Lie-rigid and $H^{2}(L,L)=HL^{2}(L,L)$  
\end{enumerate}
\end{theor}
The Leibniz algebra $L$ itself is a Leibniz module under the action $L\otimes L\longrightarrow L$ by $(x,y)\mapsto [x,y].$ In this case the two actions are given by the Leibniz multiplication.
Let $$BL^{n}(L,M)=\{f\in Hom(L^{\otimes n},M); f=d^{n}g \ for \ some \ g\in Hom(L^{\otimes n-1},M)\}$$
and $ZL^{n}(L,M)=\{f\in Hom(L^{\otimes n},M); d^{n}f=0 \}.$
Then $BL^{n}(L,M)$ is called a space of all $n$-coboundaries, while $ZL^{n}(L,M)$ is called a space of all $n$-cocycles. The fact that $d^{n+1}d^{n}=0$, $n\geq 0$ implies that $BL^{n}(L,M)$ is a subspace of $ZL^{n}(L,M).$ The quotient space 
$$HL^{n}(L,M)=ZL^{n}(L,M)/BL^{n}(L,M)$$
is called the $nth$-cohomology space of the Leibniz module $M$ 
The main field $F$ is a Leibniz module under the actions $F\otimes L\longrightarrow F$ by $(k,x)\mapsto [k,x]=0$ and $L\otimes F\longrightarrow F$ by $(x,k)\mapsto [x,k]=0$  .
An optimal way of exhaustive study of degenerations in a set of algebras includes intensive usage of necessary criteria based on degeneration invariants. The invariants are preserved under degenerations. In the next section for the further references we collect some degeneration invariants. The proofs can be found more or less in the literature, see \cite{72}.

\section{Acknowledgments}
This study was supported by the grant 05-02-12-2188RU provided by the Universiti Putra Malaysia (UPM).

\end{document}